\documentclass [a4paper,english]{article}
\usepackage{amsthm, amsmath, amssymb, amsbsy, amsfonts}
\usepackage[cp1251]{inputenc}
\usepackage[T2A]{fontenc}
\usepackage{indentfirst,epigraph}
\usepackage[english,russian]{babel}

\def\le{\leqslant}

\def\ge{\geqslant}
\def\geq{\geqslant}
\def\phi{\varphi}

\def\kappa{\varkappa}

\newtheorem{theorem}{\bf Теорема}

\newtheorem{corollary}{\bf Следствие}
\newtheorem{lemma}{\bf Лемма}
\newtheorem{remark}{\bf Замечание}

\begin{document}

\title{Gasparyan's Inequality}
\author{A.B. Pevnyi,  S.M. Sitnik}

\date{}
\maketitle

\selectlanguage{english}
\begin{abstract}

We consider a new and simpler proof of an inequality of A.S. Gasparyan, which was originally derived in terms of
complex algebraical objects --- multidimensional hyperdeterminants. Our proof is much simpler and use only standard technics such as mean inequalities. The main theorem on Gini means is essentially used. Some corollaries and problems are considered.
\end{abstract}

\newpage

\selectlanguage{russian}

\begin{center}
{\Large \bf
ОБ ОДНОМ НЕРАВЕНСТВЕ А.~С.~ГАСПАРЯНА}
\end{center}

\begin{abstract}
Предлагается новое  доказательство одного неравенства А.~С. Гаспаряна, полученного в оригинале
как следствие достаточно сложного  результата, сформулированного в терминах
некоторых нестандартных алгебраических объектов---многомерных гипердетерминантов. Мы даем упрощенное
доказательство на основе  стандартных  методов анализа с использованием средних величин и
неравенств для них. Существенно используется теорема сравнения для специального
класса средних Джини. В качестве следствия уточнены неравенства для комбинаций степенных средних.
\end{abstract}
\vskip10mm

\noindent
Гаспарян
\cite{Gas1}--\cite{Gas3}
получил несколько достаточно
простых по форме неравенств
как
следствие сложно формулируемых результатов в теории
гиперопределителей.
Приведем
одно из этих неравенств.

Пусть даны числа $m\in{\mathbb{N}}, n\in\mathbb{N}$
и
набор вещественных чисел $(x_1,\dots,x_n)$.
Введем величину
\begin{gather}
\label{f1}
G_{2m}(x_1,\dots,x_n)=\sum\limits_{r=0}^{2m}(-1)^r
\left(\sum\limits_{k=1}^n x_k^r\right)
\left(\sum\limits_{k=1}^n x_k^{2m-r}\right).
\end{gather}

\begin{theorem}[Гаспарян]
 \label{theorem_gasparyan}
Величина
$G_{2m}(x_1,\dots,x_n)$ неотрицательна,
т.е.
$G_{2m}\geq 0$.
\end{theorem}

В настоящей работе получено более простое
доказательство
теоремы \ref{theorem_gasparyan}
с использованием оценок
для средних значений, в том числе средних Джини,
а также уточнения этого неравенства.

Рассмотрим введенную  Гаспаряном степенную форму
(\ref{f1}), которая является однородным многочленом
степени $2m$ от $n$ переменных $x=(x_1,\dots,x_n)$.
Покажем, что условие неотрицательности можно усилить,
показав, что эта форма на самом деле отделена от нуля
простейшими положительными формами той
же степени $2m$.

\begin{theorem}
\label{theorem_02}
Для любого $x\in \mathbb{R}^n$ справедливы неравенства
\begin{align}
&
G_{2m}(x)\ge \left(\sum\limits_{i=1}^n  x_i^m\right)^2,
\quad
\text{если $m$~--- четное число},
\label{f2}
\\[3pt]&
G_{2m}(x)\ge \left(\sum\limits_{i=1}^n  x_i^{m-1}\right)
\left(\sum\limits_{i=1}^n  x_i^{m+1}\right),
\quad
\text{если $m$~--- нечетное число}.
\label{f3}
\end{align}
Если все $x_i$ равны между собой, то оба неравенства
\eqref{f2} и \eqref{f3} обращаются в равенство.
Обратно, если хоть одно из неравенств
\eqref{f2} и \eqref{f3}
обращается в равенство для положительных
$x_1,x_2,\dots,x_n$, то $x_1=x_2= \ldots =x_n$.
\end{theorem}

\begin{corollary} Форма из \eqref{f1} положительна,
т.е.
$G_{2m}(x)>0$ для всех $x\in \mathbb{R}^n$,
$x\not=0$.
\end{corollary}

\begin{remark}
Для правых частей в (\ref{f2}) и (\ref{f3})
справедливо неравенство
$$ \left(\sum\limits_{i=1}^n  x_i^m\right)^2 \le
\left(\sum\limits_{i=1}^n  x_i^{m-1}\right)
\left(\sum\limits_{i=1}^n  x_i^{m+1}\right) $$
для нечетного $m$. Для доказательства достаточно представить
$x_i^m$ в виде $x_i^{(m-1)/2}x_i^{(m+1)/2}$ и применить
неравенство Коши~--- Буняковского.
Поэтому оценка (\ref{f2}) выполняется для всех
значений $m$, но для нечетных значений она является
менее точной, чем (\ref{f3}).
\end{remark}

Для доказательства теоремы \ref{theorem_02} потребуются некоторые
результаты о степенных средних.

Пусть даны числа $r>0$,
$r\in \mathbb{R}$,
$m\in{\mathbb{N}}$,
$n\in\mathbb{N}$,
и набор неотрицательных чисел
$(x_1,\dots,x_n)$,
$x_k\ge 0$,
$0\le k \le n$.
Рассмотрим среднее степенное для этого набора чисел
$$
M_r(x_1,\dots,x_n)=M_r=
\left(\frac{x_1^r+\dots+x_n^r}{n}\right)^\frac{1}{r},
$$
причем тогда
$$
M_r^r=\frac{x_1^r+\ldots+x_n^r}{n},
$$
а также сумму степеней
$$
S_r(x_1,\dots,x_n)=S_r=x_1^r+\ldots+x_n^r,
$$
связанную со степенным средним соответствующего порядка очевидными формулами
$$
M_r=\left(\frac{S_r}{n}\right)^\frac{1}{r},\quad
 M_r^r=\frac{S_r}{n},\quad  S_r=n M_r^r.
$$

Свойства средних величин и степенных сумм подробно
изложены в
\cite{HLP}--\cite{MPF} и
\cite{G}--\cite{B}. Из
работ последнего времени
 отметим
\cite{Hin} и
обзоры второго автора \cite{Sit1,Sit2}.

Установим два вспомогательных утверждения.
Пусть $a_1,\dots,a_n$~--- положительные числа.
Введем величину
$$
F(r)=n^2 M_r^r M_{2m-r}^{2m-r}=S_r S_{2m-r}=
\sum\limits_{i=1}^n a_i^r \sum\limits_{j=1}^n a_j^{2m-r}.
$$

\begin{lemma}
\label{lemma_01}
Функция $F(r)$ выпукла на $(-\infty, \infty)$, не возрастает
на промежутке $0 \le r \le m$ и не убывает
на промежутке $m \le r \le 2m$. Если $F(r)=F(s)$
при некоторых $0 \le r<s \le m$, то необходимо $a_1=a_2=
\ldots
=a_n$.
\end{lemma}

\begin{proof}[\bf Доказательство]
Вычислим вторую производную
$$
F''(r)=\sum\limits_{i,j=1}^n a_i^r a_j^{2m-r} (\ln a_i - \ln a_j)^2 \ge 0.
$$
Выпуклость установлена.
Докажем неравенство $F(r)\ge F(s)$ при $0 \le r<s \le m$. Оно эквивалентно
$$
\sum\limits_{i=1}^n a_i^r
\sum\limits_{i=1}^n a_i^{2m-r}
\ge \sum\limits_{i=1}^n a_i^s \sum\limits_{i=1}^n a_i^{2m-s}.
$$
Перепишем последнее неравенство в виде отношений
$$
\frac{\sum\limits_k a_k^{2m-r}}
{\sum\limits_k a_k^{2m-s}} \ge \frac{\sum\limits_k a_k^s}{\sum\limits_k a_k^r}
$$
и возведем в положительную степень
${1}/({s-r})=
{1}/[{(2m-r)-(2m-s)}]$:
$$
\left(\frac{\sum\limits_k a_k^{2m-r}}
{\sum\limits_k a_k^{2m-s}}\right)^{{1}/[{(2m-r)-(2m-s)}]}
\ge \left(\frac{\sum\limits_k a_k^s}
{\sum\limits_k a_k^r}\right)^{{1}/({s-r})}.
$$
В результате мы получили неравенство между средними
Джини
(см. \cite{MPF}--\cite{B}). В обозначениях из
\cite{Sit1,Sit2} оно записывается так:
\begin{equation}
\label{f4}
Gi^{2m-r,2m-s} \ge Gi^{s,r}.
\end{equation}
Джини ввел и изучил средние, названные
впоследствии его именем, в \cite{G}, которые затем многократно переоткрывались.
Условия справедливости  неравенства (\ref{f4}) приведены
в \cite[с. 249]{B}.
В нашем случае
они сводятся к одному неравенству $s+r\le 2m$,
которое выполнено по предположению.

Если $F(r)=F(s)$, то, как установлено в \cite{B},
выполняются равенства $a_1=a_2=\ldots
=a_n$.

График функции $F(r)$ симметричен
относительно точки $r=m$ ввиду очевидного свойства $F(2m-r)=F(r)$.
Лемма  доказана.
\end{proof}

Покажем, что частично результаты леммы
\ref{lemma_01} можно перенести на наборы чисел любого знака.

\begin{lemma}
\label{lemma_02}
Пусть $x_1,\dots,x_n \in \mathbb{R}$.
Тогда при любом четном $r$ из промежутка $[0,m-1]$ справедливо неравенство
\begin{equation}
\label{f5}
F(r) \ge F(r+1).
\end{equation}
\end{lemma}

\begin{proof}[\bf Доказательство]
Если
$x_i>0$
для всех $i$, то неравенство (\ref{f5}) следует из
леммы \ref{lemma_01}. Пусть
$x_i$ любого знака и хоть одно $x_i$ отлично от нуля.
Тогда
\begin{eqnarray*}
F(r+1)=\sum\limits_{i=1}^n x_i^{r+1} \sum\limits_{j=1}^n x_j^{2m-r-1}
\le \sum\limits_{i=1}^n {\lvert x_i \rvert}^{r+1}
\sum\limits_{j=1}^n {\lvert x_j \rvert}^{2m-r-1} \le \\
\le \sum\limits_{i=1}^n {\lvert x_i \rvert}^{r}
\sum\limits_{j=1}^n {\lvert x_j \rvert}^{2m-r} =F(r).
\end{eqnarray*}
Здесь мы применили лемму \ref{lemma_01} к числам
$\lvert x_1 \rvert,\dots, \lvert x_n \rvert$, при
этом, если среди этих чисел есть нули, то их надо
отбросить и применить лемму \ref{lemma_01} к оставшимся числам.
Кроме того, было использовано, что $r$---четное число,
 так как тогда ${\lvert x_i \rvert}^r=x_i^r$ и
${\lvert x_i \rvert}^{2m-r}=x_i^{2m-r}$.
Лемма  доказана.
\end{proof}

Теперь можно закончить доказательство теоремы \ref{theorem_02}.
При $m$ четном получаем
$$
G_{2m}(x)=2\sum\limits_{r=0}^{m-2} (F(r)-F(r+1)) + F(m).
$$
По лемме \ref{lemma_02}
$F(r)-F(r+1)\ge 0$, поэтому $G_{2m}(x)\ge F(m)$, что совпадает с (\ref{f2}).

При $m$ нечетном получаем
$$
G_{2m}(x)=2\sum\limits_{r=0}^{m-3} (F(r)-F(r+1)) + 2F(m-1)-F(m).
$$
По лемме \ref{lemma_02}
$F(r)-F(r+1)\ge 0$ при $r=0,2,\cdots,m-1$,  поэтому
$G_{2m}(x)\ge F(m-1)$, что совпадает с (\ref{f3}).

Если неравенства (\ref{f2}) или (\ref{f3}) обращаются
в равенства, то будут выполняться соотношения вида
$F(r)-F(r+1)=0$. А тогда в случае положительных
$x_i$ по лемме \ref{lemma_01} все числа $x_i$ равны между собой.
 Это завершает доказательство теоремы \ref{theorem_02}.

Следует отметить крайние случаи значений параметра
$r$ в лемме \ref{lemma_01}. Так, при $r=0$  неравенство
$F(0)\ge F(1)$ сводится к такому:
$$
M_{2m}^{2m} \ge M_1 M_{2m-1}^{2m-1}.
$$
Это неравенство очевидно, так как сводится к
перемножению двух стандартных неравенств,
выражающих монотонность среднего степенного по
параметру (см., например, \cite{HLP,BB})
$$
M_1 \le M_{2m},\quad  M_{2m-1} \le M_{2m}.
$$

Другой крайний случай леммы \ref{lemma_01} получается из
сформулированного в ней условия минимума функции
$F(r)$ при $r=m$ и сводится к неравенству
$$
(M_{m}^{m})^2 \le M_r^r M_{2m-r}^{2m-r}, \quad
0 \le r \le m,
$$
которое является частным случаем известного неравенства
 Ляпунова (см., например, \cite{HLP,BB}),
 выражающего логарифмическую выпуклость среднего
степенного по параметру.

\vskip 20mm

АВТОРЫ:\\
\\
\\
Ситник С.М.\\
Воронежский институт МВД, Воронеж, Россия.\\
Sitnik S.M., Voronezh Institute of the Ministry of Internal Affairs of Russia.\\
mathsms@yandex.ru, pochtasms@gmail.com
\\
\\
Певный А.Б.\\
Сыктывкарский госуниверситет, Сыктывкар, Россия.\\
Pevnyi A.B., Syktyvkar State University, Syktyvkar, Russia.\\
pevnyi@syktsu.ru


\begin{thebibliography}{944}


\vskip4pt
\bibitem{Gas1}
А.С. Гаспарян, ``Аналог формулы Бине--Коши для многомерных матриц''\ ,
\textit{ДАН СССР}
\textbf{273},
No. 2, (1983).

\vskip4pt\bibitem{Gas2}
А. С. Гаспарян, ``О некоторых приложениях
 многомерных матриц''
\textit{Сообщ.
прикл. мат., ВЦ АН СССР} (1983).

\vskip4pt\bibitem{Gas3}
А. С. Гаспарян, ``Гиперопределители и
обобщенные неравенства Чебышева'',
В:
\textit{Тез. докл. межд. конф.
``Математические идеи П.Л. Чебышева
и их приложения к современным
проблемам естествознания''}
с. 32--34,
Обнинск (2002).


\vskip4pt\bibitem{HLP}
Г.Г. Харди, Дж.Е. Литтлвуд, Г. Полиа,
\textit{Неравенства}, ИЛ, М. (1948).

\vskip4pt\bibitem{BB}
Э. Беккенбах, Р. Беллман,
\textit{Неравенства}, Мир, М. (1965).

\vskip4pt\bibitem{MV}
D. S. Mitrinovi\'{c} (in cooperation with P. M. Vasi\'{c}),
\textit{Analytic Inequalities}, Springer, Berlin etc.  (1970).

\vskip4pt\bibitem{MPF}
D. S. Mitrinovi\'{c},
J. E. Pe\v{c}ari\'{c},
A. M. Fink,
\textit{Classical and New Inequalities in Analysis}, Kluwer (1993).

\vskip4pt\bibitem{G}
К. Джини,
\textit{Средние величины}, Статистика, М. (1970).

\vskip4pt\bibitem{BMV}
P. S. Bullen, D. S. Mitrinovi\'{c}, P. M. Vasi\'{c},
\textit{Means and Their Inequalities},
D.~Reidel Publ., Dordrecht (1988).

\vskip4pt\bibitem{B}
P. S. Bullen,
\textit{Handbook of Means and Their Inequalities}, Kluwer (2003).

\vskip4pt\bibitem{Hin}
M. Hajja, P. S. Bullen, J. Matkowski, E. Neuman, S. Simic (Eds.),
\textit{Means and Their Inequalities}.
Special issue of Int. J. Math.
Math. Sci. (2013).

\vskip4pt\bibitem{Sit1}
С. М. Ситник, ``Уточнения и обобщения
классических неравенств'',
В:
\textit{Итоги науки}, с. 221--266,
Владикавказ (2009).

\vskip4pt\bibitem{Sit2}
S. M. Sitnik,  ``Generalized Young and
Cauchy--Bunyakowsky Inequalities with
  Applications: A Survey'', \textit{arXiv:1012.3864}.




\end{thebibliography}
\end{document}